\documentclass{gtart}
\def\ifplaintex{\expandafter\ifx\csname documentclass\endcsname\relax}

\def\gtp{{\mathsurround=0pt\it $\cal G\mskip-2mu$eometry \&\ 
$\cal T\!\!$opology $\cal P\!$ublications}}  % GT publications

\def\recd{{\small Received:\qua\receiveddate\ifx\reviseddate\relax
\else\qquad Revised:\qua\reviseddate\fi\par}} 

\def\lognumber#1{\def\thelognumber{#1}}
\def\volumenumber#1{\def\thevolumenumber{#1}}
\def\volumeyear#1{\def\thevolumeyear{#1}}
\def\papernumber#1{\def\thepapernumber{#1}}
\def\pagenumbers#1#2{\def\startpage{#1}\def\finishpage{#2}}
\def\published#1{\def\publishdate{#1}}

\def\received#1{\def\receiveddate{#1}}
\def\revised#1{\def\reviseddate{#1}}
\def\accepted#1{\def\accepteddate{#1}}

\long\def\asciiabstract#1{\long\def\theasciiabstract{#1}}

\let\\\par\let\thelognumber\relax\let\thevolumenumber\relax
\let\thepapernumber\relax\let\thevolumeyear\relax\let\startpage\relax
\let\finishpage\relax\let\publishdate\relax\let\receiveddate\relax
\let\reviseddate\relax\let\accepteddate\relax\let\theasciititle\relax
\let\theasciiauthors\relax
\let\theasciiabstract\relax

\let\theasciiemail\relax

\ifplaintex
\font\logobig=cmssbx10 scaled 3836
\font\logomed=cmssbx10 scaled 2557
\else
\font\logobig=cmssbx10 scaled 4200
\font\logomed=cmssbx10 scaled 2800
\fi

\long\def\makeagttitle{   %%% start of definition of \makeagttitle
\count0=\startpage
\agt\hfill      %   Journal title (top left) 
\hbox to 45truept{\vbox to 0pt{\vglue -13truept{\logomed A\kern -.37em{\logobig 
T}\kern -.38em G}\vss}\hss}
\break
{\small Volume \thevolumenumber\ (\thevolumeyear)
\startpage--\finishpage\nl
Published: \publishdate}

\vglue .25truein

{\parskip=0pt\leftskip 0pt plus
1fil\def\\{\par\smallskip}{\Large\bf\thetitle}\par\medskip} \vglue
0.05truein

{\parskip=0pt\leftskip 0pt plus 1fil\def\\{\par}{\sc\theauthors}
\par\medskip}%
 
\vglue 0.03truein 

{\small\leftskip 25truept\rightskip 25truept{\bf Abstract}\stdspace\theabstract

{\bf AMS Classification}\stdspace\theprimaryclass
\ifx\thesecondaryclass\relax\else; \thesecondaryclass\fi\par
{\bf Keywords}\stdspace \thekeywords\par}\vglue 7truept

}   %%%% end of definition of \makeagttitle

\ifplaintex
\font\phead=cmsl9 scaled 950
\font\pnum=cmbx10 scaled 913
\font\pfoot=cmsl9 scaled 950
\headline{\vbox to 0pt{\vskip -4.5mm\line{\small\phead\ifnum
\count0=\startpage ISSN 1472-2739 (on-line) 1472-2747 (printed)
\hfill {\pnum\folio}\else\ifodd\count0\def\\{ }% 
\ifx\theshorttitle\relax\thetitle\else\theshorttitle\fi\hfill{\pnum\folio}
\else\def\\{ and }{\pnum\folio}\hfill\ifx\theshortauthors\relax\theauthors
\else\theshortauthors\fi\fi\fi}\vss}}
\footline{\vbox to 0pt{\vglue 0mm\line{\small\pfoot\ifnum\count0=\startpage
\copyright\ \gtp\hfill\else
\agt, Volume \thevolumenumber\ (\thevolumeyear)\hfill\fi}\vss}}
\else
\font=cmsl9 scaled 1050
\font\lnum=cmbx10 
\font=cmsl9 scaled 1050
\makeatletter
\def\@oddhead{{\small\lhead\ifnum\count0=\startpage ISSN 1472-2739 
(on-line) 1472-2747 (printed)\hfill {\lnum\number\count0}\else\ifodd\count0
\def\\{ }\ifx\theshorttitle\relax \thetitle \else\theshorttitle\fi\hfill
{\lnum\number\count0}\else\def\\{ and }{\lnum\number\count0}
\hfill\ifx\theshortauthors\relax 
\theauthors\else\theshortauthors\fi\fi\fi}}\def\@evenhead{\@oddhead}
\def\@oddfoot{\small\lfoot\ifnum\count0=\startpage\copyright\ \gtp\hfill\else
\agt, Volume \thevolumenumber\ (\thevolumeyear)\hfill\fi}
\def\@evenfoot{\@oddfoot}
\makeatother
\fi
\let\maketitlepage\makeagttitle
\let\makeshorttitle\maketitlepage
\let\maketitle\maketitlepage

\newwrite\gtoutfile
\long\gdef\makeheadfile{  %%% start of definition of \makeheadfile
{\def\\{, }\def\s{ }
\immediate\openout\gtoutfile head.xxx
\immediate\write\gtoutfile{To: math@arxiv.org}
\immediate\write\gtoutfile{Subject: put OR rep NNNNN:ppppp}
\immediate\write\gtoutfile{--text follows this line--}
\immediate\write\gtoutfile{Proxy-for: \ifx\theasciiauthors\relax
\theauthors\else\theasciiauthors\fi\s<\ifx\theasciiemail\relax\theemail\else\theasciiemail\fi>}
\immediate\write\gtoutfile{\noexpand\\}
\immediate\write\gtoutfile{Authors: \ifx\theasciiauthors\relax
\theauthors\else\theasciiauthors\fi}
{\def\\{ }\immediate\write\gtoutfile{Title: \ifx\theasciititle\relax
\thetitle\else\theasciititle\fi}}
\immediate\write\gtoutfile{Subj-class: GT or SG, GR etc}
\immediate\write\gtoutfile{MSC-class: \theprimaryclass\ifx\thesecondaryclass\relax\else, \thesecondaryclass\fi}
\immediate\write\gtoutfile{Journal-ref: Algebr. Geom. Topol. \thevolumenumber\s
(\thevolumeyear) \startpage-\finishpage}
\immediate\write\gtoutfile{Comments: Published by Algebraic and
Geometric Topology at}
\immediate\write\gtoutfile{\s\s\s  http://www.maths.warwick.ac.uk/agt/AGTVol\thevolumenumber/agt-\thevolumenumber-\thepapernumber.abs.html}
\immediate\write\gtoutfile{\noexpand\\}
\immediate\write\gtoutfile{}
\ifx\theasciiabstract\relax
\immediate\write\gtoutfile{\theabstract}\else
\immediate\write\gtoutfile{\theasciiabstract}\fi
\immediate\write\gtoutfile{}
\immediate\write\gtoutfile{\noexpand\\}
\immediate\write\gtoutfile{}
\immediate\closeout\gtoutfile}}  %%% end of definition of \makeheadfile

\def\maketitlepage{\makeagttitle\makeheadfile}
\let\makeshorttitle\maketitlepage
\let\maketitle\maketitlepage

\def\ifplaintex{\expandafter\ifx\csname documentclass\endcsname\relax}

\def\gtp{{\mathsurround=0pt\it $\cal G\mskip-2mu$eometry \&\ 
$\cal T\!\!$opology $\cal P\!$ublications}}  % GT publications

\def\recd{{\small Received:\qua\receiveddate\ifx\reviseddate\relax
\else\qquad Revised:\qua\reviseddate\fi\par}} 

\def\lognumber#1{\def\thelognumber{#1}}
\def\volumenumber#1{\def\thevolumenumber{#1}}
\def\volumeyear#1{\def\thevolumeyear{#1}}
\def\papernumber#1{\def\thepapernumber{#1}}
\def\pagenumbers#1#2{\def\startpage{#1}\def\finishpage{#2}}
\def\published#1{\def\publishdate{#1}}

\def\received#1{\def\receiveddate{#1}}
\def\revised#1{\def\reviseddate{#1}}
\def\accepted#1{\def\accepteddate{#1}}

\long\def\asciiabstract#1{\long\def\theasciiabstract{#1}}

\let\\\par\let\thelognumber\relax\let\thevolumenumber\relax
\let\thepapernumber\relax\let\thevolumeyear\relax\let\startpage\relax
\let\finishpage\relax\let\publishdate\relax\let\receiveddate\relax
\let\reviseddate\relax\let\accepteddate\relax\let\theasciititle\relax
\let\theasciiauthors\relax
\let\theasciiabstract\relax

\let\theasciiemail\relax

\ifplaintex
\font\logobig=cmssbx10 scaled 3836
\font\logomed=cmssbx10 scaled 2557
\else
\font\logobig=cmssbx10 scaled 4200
\font\logomed=cmssbx10 scaled 2800
\fi

\long\def\makeagttitle{   %%% start of definition of \makeagttitle
\count0=\startpage
\agt\hfill      %   Journal title (top left) 
\hbox to 45truept{\vbox to 0pt{\vglue -13truept{\logomed A\kern -.37em{\logobig 
T}\kern -.38em G}\vss}\hss}
\break
{\small Volume \thevolumenumber\ (\thevolumeyear)
\startpage--\finishpage\nl
Published: \publishdate}

\vglue .25truein

{\parskip=0pt\leftskip 0pt plus
1fil\def\\{\par\smallskip}{\Large\bf\thetitle}\par\medskip} \vglue
0.05truein

{\parskip=0pt\leftskip 0pt plus 1fil\def\\{\par}{\sc\theauthors}
\par\medskip}%
 
\vglue 0.03truein 

{\small\leftskip 25truept\rightskip 25truept{\bf Abstract}\stdspace\theabstract

{\bf AMS Classification}\stdspace\theprimaryclass
\ifx\thesecondaryclass\relax\else; \thesecondaryclass\fi\par
{\bf Keywords}\stdspace \thekeywords\par}\vglue 7truept

}   %%%% end of definition of \makeagttitle

\ifplaintex
\font\phead=cmsl9 scaled 950
\font\pnum=cmbx10 scaled 913
\font\pfoot=cmsl9 scaled 950
\headline{\vbox to 0pt{\vskip -4.5mm\line{\small\phead\ifnum
\count0=\startpage ISSN 1472-2739 (on-line) 1472-2747 (printed)
\hfill {\pnum\folio}\else\ifodd\count0\def\\{ }% 
\ifx\theshorttitle\relax\thetitle\else\theshorttitle\fi\hfill{\pnum\folio}
\else\def\\{ and }{\pnum\folio}\hfill\ifx\theshortauthors\relax\theauthors
\else\theshortauthors\fi\fi\fi}\vss}}
\footline{\vbox to 0pt{\vglue 0mm\line{\small\pfoot\ifnum\count0=\startpage
\copyright\ \gtp\hfill\else
\agt, Volume \thevolumenumber\ (\thevolumeyear)\hfill\fi}\vss}}
\else
\font=cmsl9 scaled 1050
\font\lnum=cmbx10 
\font=cmsl9 scaled 1050
\makeatletter
\def\@oddhead{{\small\lhead\ifnum\count0=\startpage ISSN 1472-2739 
(on-line) 1472-2747 (printed)\hfill {\lnum\number\count0}\else\ifodd\count0
\def\\{ }\ifx\theshorttitle\relax \thetitle \else\theshorttitle\fi\hfill
{\lnum\number\count0}\else\def\\{ and }{\lnum\number\count0}
\hfill\ifx\theshortauthors\relax 
\theauthors\else\theshortauthors\fi\fi\fi}}\def\@evenhead{\@oddhead}
\def\@oddfoot{\small\lfoot\ifnum\count0=\startpage\copyright\ \gtp\hfill\else
\agt, Volume \thevolumenumber\ (\thevolumeyear)\hfill\fi}
\def\@evenfoot{\@oddfoot}
\makeatother
\fi
\let\maketitlepage\makeagttitle
\let\makeshorttitle\maketitlepage
\let\maketitle\maketitlepage

\newwrite\gtoutfile
\long\gdef\makeheadfile{  %%% start of definition of \makeheadfile
{\def\\{, }\def\s{ }
\immediate\openout\gtoutfile head.xxx
\immediate\write\gtoutfile{To: math@arxiv.org}
\immediate\write\gtoutfile{Subject: put OR rep NNNNN:ppppp}
\immediate\write\gtoutfile{--text follows this line--}
\immediate\write\gtoutfile{Proxy-for: \ifx\theasciiauthors\relax
\theauthors\else\theasciiauthors\fi\s<\ifx\theasciiemail\relax\theemail\else\theasciiemail\fi>}
\immediate\write\gtoutfile{\noexpand\\}
\immediate\write\gtoutfile{Authors: \ifx\theasciiauthors\relax
\theauthors\else\theasciiauthors\fi}
{\def\\{ }\immediate\write\gtoutfile{Title: \ifx\theasciititle\relax
\thetitle\else\theasciititle\fi}}
\immediate\write\gtoutfile{Subj-class: GT or SG, GR etc}
\immediate\write\gtoutfile{MSC-class: \theprimaryclass\ifx\thesecondaryclass\relax\else, \thesecondaryclass\fi}
\immediate\write\gtoutfile{Journal-ref: Algebr. Geom. Topol. \thevolumenumber\s
(\thevolumeyear) \startpage-\finishpage}
\immediate\write\gtoutfile{Comments: Published by Algebraic and
Geometric Topology at}
\immediate\write\gtoutfile{\s\s\s  http://www.maths.warwick.ac.uk/agt/AGTVol\thevolumenumber/agt-\thevolumenumber-\thepapernumber.abs.html}
\immediate\write\gtoutfile{\noexpand\\}
\immediate\write\gtoutfile{}
\ifx\theasciiabstract\relax
\immediate\write\gtoutfile{\theabstract}\else
\immediate\write\gtoutfile{\theasciiabstract}\fi
\immediate\write\gtoutfile{}
\immediate\write\gtoutfile{\noexpand\\}
\immediate\write\gtoutfile{}
\immediate\closeout\gtoutfile}}  %%% end of definition of \makeheadfile

\def\maketitlepage{\makeagttitle\makeheadfile}
\let\makeshorttitle\maketitlepage
\let\maketitle\maketitlepage

\def\ifplaintex{\expandafter\ifx\csname documentclass\endcsname\relax}

\def\gtp{{\mathsurround=0pt\it $\cal G\mskip-2mu$eometry \&\ 
$\cal T\!\!$opology $\cal P\!$ublications}}  % GT publications

\def\recd{{\small Received:\qua\receiveddate\ifx\reviseddate\relax
\else\qquad Revised:\qua\reviseddate\fi\par}} 

\def\lognumber#1{\def\thelognumber{#1}}
\def\volumenumber#1{\def\thevolumenumber{#1}}
\def\volumeyear#1{\def\thevolumeyear{#1}}
\def\papernumber#1{\def\thepapernumber{#1}}
\def\pagenumbers#1#2{\def\startpage{#1}\def\finishpage{#2}}
\def\published#1{\def\publishdate{#1}}

\def\received#1{\def\receiveddate{#1}}
\def\revised#1{\def\reviseddate{#1}}
\def\accepted#1{\def\accepteddate{#1}}

\long\def\asciiabstract#1{\long\def\theasciiabstract{#1}}

\let\\\par\let\thelognumber\relax\let\thevolumenumber\relax
\let\thepapernumber\relax\let\thevolumeyear\relax\let\startpage\relax
\let\finishpage\relax\let\publishdate\relax\let\receiveddate\relax
\let\reviseddate\relax\let\accepteddate\relax\let\theasciititle\relax
\let\theasciiauthors\relax
\let\theasciiabstract\relax

\let\theasciiemail\relax

\ifplaintex
\font\logobig=cmssbx10 scaled 3836
\font\logomed=cmssbx10 scaled 2557
\else
\font\logobig=cmssbx10 scaled 4200
\font\logomed=cmssbx10 scaled 2800
\fi

\long\def\makeagttitle{   %%% start of definition of \makeagttitle
\count0=\startpage
\agt\hfill      %   Journal title (top left) 
\hbox to 45truept{\vbox to 0pt{\vglue -13truept{\logomed A\kern -.37em{\logobig 
T}\kern -.38em G}\vss}\hss}
\break
{\small Volume \thevolumenumber\ (\thevolumeyear)
\startpage--\finishpage\nl
Published: \publishdate}

\vglue .25truein

{\parskip=0pt\leftskip 0pt plus
1fil\def\\{\par\smallskip}{\Large\bf\thetitle}\par\medskip} \vglue
0.05truein

{\parskip=0pt\leftskip 0pt plus 1fil\def\\{\par}{\sc\theauthors}
\par\medskip}%
 
\vglue 0.03truein 

{\small\leftskip 25truept\rightskip 25truept{\bf Abstract}\stdspace\theabstract

{\bf AMS Classification}\stdspace\theprimaryclass
\ifx\thesecondaryclass\relax\else; \thesecondaryclass\fi\par
{\bf Keywords}\stdspace \thekeywords\par}\vglue 7truept

}   %%%% end of definition of \makeagttitle

\ifplaintex
\font\phead=cmsl9 scaled 950
\font\pnum=cmbx10 scaled 913
\font\pfoot=cmsl9 scaled 950
\headline{\vbox to 0pt{\vskip -4.5mm\line{\small\phead\ifnum
\count0=\startpage ISSN 1472-2739 (on-line) 1472-2747 (printed)
\hfill {\pnum\folio}\else\ifodd\count0\def\\{ }% 
\ifx\theshorttitle\relax\thetitle\else\theshorttitle\fi\hfill{\pnum\folio}
\else\def\\{ and }{\pnum\folio}\hfill\ifx\theshortauthors\relax\theauthors
\else\theshortauthors\fi\fi\fi}\vss}}
\footline{\vbox to 0pt{\vglue 0mm\line{\small\pfoot\ifnum\count0=\startpage
\copyright\ \gtp\hfill\else
\agt, Volume \thevolumenumber\ (\thevolumeyear)\hfill\fi}\vss}}
\else
\font=cmsl9 scaled 1050
\font\lnum=cmbx10 
\font=cmsl9 scaled 1050
\makeatletter
\def\@oddhead{{\small\lhead\ifnum\count0=\startpage ISSN 1472-2739 
(on-line) 1472-2747 (printed)\hfill {\lnum\number\count0}\else\ifodd\count0
\def\\{ }\ifx\theshorttitle\relax \thetitle \else\theshorttitle\fi\hfill
{\lnum\number\count0}\else\def\\{ and }{\lnum\number\count0}
\hfill\ifx\theshortauthors\relax 
\theauthors\else\theshortauthors\fi\fi\fi}}\def\@evenhead{\@oddhead}
\def\@oddfoot{\small\lfoot\ifnum\count0=\startpage\copyright\ \gtp\hfill\else
\agt, Volume \thevolumenumber\ (\thevolumeyear)\hfill\fi}
\def\@evenfoot{\@oddfoot}
\makeatother
\fi
\let\maketitlepage\makeagttitle
\let\makeshorttitle\maketitlepage
\let\maketitle\maketitlepage

\newwrite\gtoutfile
\long\gdef\makeheadfile{  %%% start of definition of \makeheadfile
{\def\\{, }\def\s{ }
\immediate\openout\gtoutfile head.xxx
\immediate\write\gtoutfile{To: math@arxiv.org}
\immediate\write\gtoutfile{Subject: put OR rep NNNNN:ppppp}
\immediate\write\gtoutfile{--text follows this line--}
\immediate\write\gtoutfile{Proxy-for: \ifx\theasciiauthors\relax
\theauthors\else\theasciiauthors\fi\s<\ifx\theasciiemail\relax\theemail\else\theasciiemail\fi>}
\immediate\write\gtoutfile{\noexpand\\}
\immediate\write\gtoutfile{Authors: \ifx\theasciiauthors\relax
\theauthors\else\theasciiauthors\fi}
{\def\\{ }\immediate\write\gtoutfile{Title: \ifx\theasciititle\relax
\thetitle\else\theasciititle\fi}}
\immediate\write\gtoutfile{Subj-class: GT or SG, GR etc}
\immediate\write\gtoutfile{MSC-class: \theprimaryclass\ifx\thesecondaryclass\relax\else, \thesecondaryclass\fi}
\immediate\write\gtoutfile{Journal-ref: Algebr. Geom. Topol. \thevolumenumber\s
(\thevolumeyear) \startpage-\finishpage}
\immediate\write\gtoutfile{Comments: Published by Algebraic and
Geometric Topology at}
\immediate\write\gtoutfile{\s\s\s  http://www.maths.warwick.ac.uk/agt/AGTVol\thevolumenumber/agt-\thevolumenumber-\thepapernumber.abs.html}
\immediate\write\gtoutfile{\noexpand\\}
\immediate\write\gtoutfile{}
\ifx\theasciiabstract\relax
\immediate\write\gtoutfile{\theabstract}\else
\immediate\write\gtoutfile{\theasciiabstract}\fi
\immediate\write\gtoutfile{}
\immediate\write\gtoutfile{\noexpand\\}
\immediate\write\gtoutfile{}
\immediate\closeout\gtoutfile}}  %%% end of definition of \makeheadfile

\def\maketitlepage{\makeagttitle\makeheadfile}
\let\makeshorttitle\maketitlepage
\let\maketitle\maketitlepage

\def\ifplaintex{\expandafter\ifx\csname documentclass\endcsname\relax}

\def\gtp{{\mathsurround=0pt\it $\cal G\mskip-2mu$eometry \&\ 
$\cal T\!\!$opology $\cal P\!$ublications}}  % GT publications

\def\recd{{\small Received:\qua\receiveddate\ifx\reviseddate\relax
\else\qquad Revised:\qua\reviseddate\fi\par}} 

\def\lognumber#1{\def\thelognumber{#1}}
\def\volumenumber#1{\def\thevolumenumber{#1}}
\def\volumeyear#1{\def\thevolumeyear{#1}}
\def\papernumber#1{\def\thepapernumber{#1}}
\def\pagenumbers#1#2{\def\startpage{#1}\def\finishpage{#2}}
\def\published#1{\def\publishdate{#1}}

\def\received#1{\def\receiveddate{#1}}
\def\revised#1{\def\reviseddate{#1}}
\def\accepted#1{\def\accepteddate{#1}}

\long\def\asciiabstract#1{\long\def\theasciiabstract{#1}}

\let\\\par\let\thelognumber\relax\let\thevolumenumber\relax
\let\thepapernumber\relax\let\thevolumeyear\relax\let\startpage\relax
\let\finishpage\relax\let\publishdate\relax\let\receiveddate\relax
\let\reviseddate\relax\let\accepteddate\relax\let\theasciititle\relax
\let\theasciiauthors\relax
\let\theasciiabstract\relax

\let\theasciiemail\relax

\ifplaintex
\font\logobig=cmssbx10 scaled 3836
\font\logomed=cmssbx10 scaled 2557
\else
\font\logobig=cmssbx10 scaled 4200
\font\logomed=cmssbx10 scaled 2800
\fi

\long\def\makeagttitle{   %%% start of definition of \makeagttitle
\count0=\startpage
\agt\hfill      %   Journal title (top left) 
\hbox to 45truept{\vbox to 0pt{\vglue -13truept{\logomed A\kern -.37em{\logobig 
T}\kern -.38em G}\vss}\hss}
\break
{\small Volume \thevolumenumber\ (\thevolumeyear)
\startpage--\finishpage\nl
Published: \publishdate}

\vglue .25truein

{\parskip=0pt\leftskip 0pt plus
1fil\def\\{\par\smallskip}{\Large\bf\thetitle}\par\medskip} \vglue
0.05truein

{\parskip=0pt\leftskip 0pt plus 1fil\def\\{\par}{\sc\theauthors}
\par\medskip}%
 
\vglue 0.03truein 

{\small\leftskip 25truept\rightskip 25truept{\bf Abstract}\stdspace\theabstract

{\bf AMS Classification}\stdspace\theprimaryclass
\ifx\thesecondaryclass\relax\else; \thesecondaryclass\fi\par
{\bf Keywords}\stdspace \thekeywords\par}\vglue 7truept

}   %%%% end of definition of \makeagttitle

\ifplaintex
\font\phead=cmsl9 scaled 950
\font\pnum=cmbx10 scaled 913
\font\pfoot=cmsl9 scaled 950
\headline{\vbox to 0pt{\vskip -4.5mm\line{\small\phead\ifnum
\count0=\startpage ISSN 1472-2739 (on-line) 1472-2747 (printed)
\hfill {\pnum\folio}\else\ifodd\count0\def\\{ }% 
\ifx\theshorttitle\relax\thetitle\else\theshorttitle\fi\hfill{\pnum\folio}
\else\def\\{ and }{\pnum\folio}\hfill\ifx\theshortauthors\relax\theauthors
\else\theshortauthors\fi\fi\fi}\vss}}
\footline{\vbox to 0pt{\vglue 0mm\line{\small\pfoot\ifnum\count0=\startpage
\copyright\ \gtp\hfill\else
\agt, Volume \thevolumenumber\ (\thevolumeyear)\hfill\fi}\vss}}
\else
\font=cmsl9 scaled 1050
\font\lnum=cmbx10 
\font=cmsl9 scaled 1050
\makeatletter
\def\@oddhead{{\small\lhead\ifnum\count0=\startpage ISSN 1472-2739 
(on-line) 1472-2747 (printed)\hfill {\lnum\number\count0}\else\ifodd\count0
\def\\{ }\ifx\theshorttitle\relax \thetitle \else\theshorttitle\fi\hfill
{\lnum\number\count0}\else\def\\{ and }{\lnum\number\count0}
\hfill\ifx\theshortauthors\relax 
\theauthors\else\theshortauthors\fi\fi\fi}}\def\@evenhead{\@oddhead}
\def\@oddfoot{\small\lfoot\ifnum\count0=\startpage\copyright\ \gtp\hfill\else
\agt, Volume \thevolumenumber\ (\thevolumeyear)\hfill\fi}
\def\@evenfoot{\@oddfoot}
\makeatother
\fi
\let\maketitlepage\makeagttitle
\let\makeshorttitle\maketitlepage
\let\maketitle\maketitlepage

\newwrite\gtoutfile
\long\gdef\makeheadfile{  %%% start of definition of \makeheadfile
{\def\\{, }\def\s{ }
\immediate\openout\gtoutfile head.xxx
\immediate\write\gtoutfile{To: math@arxiv.org}
\immediate\write\gtoutfile{Subject: put OR rep NNNNN:ppppp}
\immediate\write\gtoutfile{--text follows this line--}
\immediate\write\gtoutfile{Proxy-for: \ifx\theasciiauthors\relax
\theauthors\else\theasciiauthors\fi\s<\ifx\theasciiemail\relax\theemail\else\theasciiemail\fi>}
\immediate\write\gtoutfile{\noexpand\\}
\immediate\write\gtoutfile{Authors: \ifx\theasciiauthors\relax
\theauthors\else\theasciiauthors\fi}
{\def\\{ }\immediate\write\gtoutfile{Title: \ifx\theasciititle\relax
\thetitle\else\theasciititle\fi}}
\immediate\write\gtoutfile{Subj-class: GT or SG, GR etc}
\immediate\write\gtoutfile{MSC-class: \theprimaryclass\ifx\thesecondaryclass\relax\else, \thesecondaryclass\fi}
\immediate\write\gtoutfile{Journal-ref: Algebr. Geom. Topol. \thevolumenumber\s
(\thevolumeyear) \startpage-\finishpage}
\immediate\write\gtoutfile{Comments: Published by Algebraic and
Geometric Topology at}
\immediate\write\gtoutfile{\s\s\s  http://www.maths.warwick.ac.uk/agt/AGTVol\thevolumenumber/agt-\thevolumenumber-\thepapernumber.abs.html}
\immediate\write\gtoutfile{\noexpand\\}
\immediate\write\gtoutfile{}
\ifx\theasciiabstract\relax
\immediate\write\gtoutfile{\theabstract}\else
\immediate\write\gtoutfile{\theasciiabstract}\fi
\immediate\write\gtoutfile{}
\immediate\write\gtoutfile{\noexpand\\}
\immediate\write\gtoutfile{}
\immediate\closeout\gtoutfile}}  %%% end of definition of \makeheadfile

\def\maketitlepage{\makeagttitle\makeheadfile}
\let\makeshorttitle\maketitlepage
\let\maketitle\maketitlepage

\lognumber{29}
\volumenumber{2}
\volumeyear{2002}
\papernumber{29}
\published{9 August 2002}
\pagenumbers{649}{664}
\received{9 May 2001}
\revised{17 April 2002}
\accepted{20 June 2002}

\usepackage{amsmath,amssymb,epsf}

\renewcommand{\hat}{\widehat}
\renewcommand{\tilde}{\widetilde}

\newfont{\cyr}{wncyr10 at 11pt}

\newcommand{\A}{\mathcal A}
\newcommand{\knots}{{\bf knots}}

\newcommand{\blob}{*}
\newcommand{\Z}{\mathbb{Z}}

\newcommand{\Q}{\mathbb{Q}}
\newcommand{\ch}{\operatorname{ch}}
\newcommand{\Vect}{\operatorname{Vect}}

\newtheorem{theorem}{Theorem} 

\newtheorem{cor}[theorem]{Corollary}
\newtheorem{prop}[theorem]{Proposition} 

\newtheorem*{qcor}{Corollary \ref{LieIntegrality}}
\newtheorem*{qthm}{Theorem \ref{DiagramNonIntegrality}}

\theoremstyle{remark}

 \newcommand{\f}{{{}^f\!\!}}
\newcommand{\thetadiag}{{\vpic{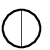}}}

\newcommand{\vpic}[1]{\vcenter{\hbox{\epsfbox{#1.eps}}}}

\title{An almost-integral universal Vassiliev\\invariant of knots} 
\author{Simon Willerton}
\email{S.Willerton@sheffield.ac.uk}
\address{Department of Pure Mathematics, University of Sheffield\\The Hicks
Building, Hounsfield Road, Sheffield, S3 7RH, UK} 
\begin{document}

\begin{abstract}
 A ``total Chern class'' invariant of knots is defined.  This is a
 universal Vassiliev invariant which is integral ``on the level of Lie
 algebras'' but it is not expressible as an integer sum of diagrams.
 The construction is motivated by similarities between the Kontsevich
 integral and the topological Chern character.
\end{abstract}
\asciiabstract{
 A `total Chern class' invariant of knots is defined.  This is a
 universal Vassiliev invariant which is integral `on the level of Lie
 algebras' but it is not expressible as an integer sum of diagrams.
 The construction is motivated by similarities between the Kontsevich
 integral and the topological Chern character.}

\maketitle

\section*{Introduction}
The Kontsevich integral is a knot invariant with two related, but
different universal properties.  The first is that it is a universal
rational Vassiliev invariant (see
Section~\ref{Section:ChernVKontsevich}), and the second is that it is
universal for perturbative quantum knot invariants coming from simple
Lie algebras (see Section~\ref{Section:LieAlgebraWeights}).  It can be
viewed as a map, $Z\colon
\Z\knots\to \A$ from the abelian group of integer
linear combinations of oriented knots to the rational algebra $\A$
which is generated by diagrams and is subject to the 1T and STU
relations (see \cite{Bar:OnVIs}).  For example,
$$Z\left(\vpic{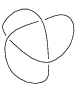} \right)=\vpic{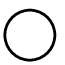}-\vpic{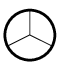}+\vpic{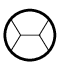}
   -\tfrac{31}{24}\vpic{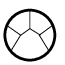}+\tfrac{5}{24}\vpic{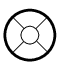}
   + \tfrac{1}{2}\vpic{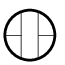} +\dots.$$
There is a lot of algebraic structure living in the algebra $\A$ (for a
survey see \cite{Willerton:KontIntAlgStructures}) some of which does
reflect topological structure of knots.  However there is no current
topological interpretation of either the Kontsevich integral or the
algebra $\A$.

The central observation of this paper is that several properties of
the Kontsevich integral are remiscent of the topological Chern
character.  Some of these are summarized in
Table~\ref{Table:ChernVKontsevich} and details are given in
Section~\ref{Section:ChernVKontsevich}.

%%%%%%%%%%%%%
%%% Table %%%
%%%%%%%%%%%%%

\begin{table}
\renewcommand{\arraystretch}{1.5}
\begin{tabular}{|c|c|}\hline
\sc Chern character&\sc Kontsevich integral\\
\hline\hline
\multicolumn{2}{|c|}{Ring map} \\\hline
$\ch\colon K(X)\to H^*(X;\Q)$&  $Z\colon \Z\knots\to \A$\\\hline
\multicolumn{2}{|c|}{Associated rational graded map is an isomorphism}\\\hline
$(K_{2i}(X)/K_{2i+1})\otimes\Q\cong H^{2i}(X;\Q)$&
$\bigl( (\Z\knots)_{i}/(\Z\knots)_{i+1}\bigr) \otimes \Q
  \cong \A_{i}$\\\hline
\multicolumn{2}{|c|}{Integrality --- when $H^{2i}(X,\Z)$ (resp.\
$\A^\Z$) is torsion-free}\\
\hline
\parbox[c][1.5\height]{15em}{
%(If $ H^{2i}(X,\Z)$ is torsion-free)\\
$x\in H^{2i}(X,\Z)\subset H^{2i}(X,\Q)\Leftrightarrow \\
\hbox{\ \ } \exists
\eta\in K_{2i}(X) \text{ st } \ch(\eta)=x+\text{hot}$}&
\parbox{15.5em}{
%(If $\A^\Z$ is torsion free)\\
$D\in\A^\Z_i\subset A_i \Leftrightarrow\hfill \\
\hbox{\ \ } \exists k\in (\Z\knots)_{i} \ \text{st}\ Z(k)=D+\text{hot}$}\\
\hline
\multicolumn{2}{|c|}{Riemann-Roch theorem v.\ cabling formula}\\
\hline
\parbox[c][1.5\height]{15.5em}{$\ch \bigl(f_K(\eta)\bigr)\cup \hat A(Y)^{-1}=\hfill\\
\hbox{\ } f_H\Bigl(\ch(\eta)\cup \hat  A(X)^{-1}\cup \exp(c_1(f)/2)\Bigr)$}&
\parbox{15em}{$\f Z\left(\Psi^{m,p}(k)\right).\Omega=\hfill\\
\hbox{\qquad}\psi^m\left(\f Z(k).\Omega.
           \exp\left(\tfrac{p}{2m}\thetadiag\right)\right)$}\\
\hline
\end{tabular}
\caption{Analogies between the topological Chern character and the
Kontsevich integral for knots.  Details are given in
Section~\ref{Section:ChernVKontsevich}.} 
\label{Table:ChernVKontsevich}
\end{table}

%%%%%%%%%%%%

Given this observation, an immediate question is ``Is there a space
$\mathbb X$ for which the Kontsevich integral\/ {\em is\/} the Chern
character?''  A positive answer to this would certainly shed an
interesting light on the Kontsevich integral, but to expect such a
space seems a little too optimistic, for reasons given below.  One way
to test such an idea is to look at consequences of the existence of
such a space.  One consequence is the existence of a total Chern
class.  The topological total Chern class would be a group homomorphism
$K(\mathbb X)\to H^*(\mathbb X,\Z)$ obtained by summing the Chern
classes together.  Note in particular that this takes values in\/ {\em
integer\/} cohomology, and the terms in it are closer to actual
topological objects than the terms in the Chern character, although
given the Chern character the total Chern class (modulo torsion) is
easily obtained.  This fact is used in
Section~\ref{SubSection:KnotChernClass} to define a ``total Chern
class for knots'' $c\colon \Z\knots\to \A$.  The obvious question to
ask at this point is ``Is this integer valued?''  The answer to this
is yes and no, as I will now explain.

This integrality question can be asked ``on the level of simple Lie
algebras'' as is done in Section~\ref{Section:LieAlgebraWeights}.
This means looking at the total Chern class after applying a weight
system coming from an irreducible representation of a simple Lie
algebra.  In this case we do get integer values:

\begin{qcor}
 The image of the total Chern class of knots is integral as far as
  irreducible representations of simple Lie algebras can detect, i.e.\ 
  if $\rho$ is such a representation, then $w_\rho\circ
  c(\Z\knots)\subset \Z[[h]]$.
\end{qcor}

%This is proved by showing that $w_\rho\circ c$ is the same as taking
%the quantum invariant, living in $\Z[q^{\pm 1}]$, associated to the
%representation $\rho$ and then applying the strange group homomorphism
%$\Z[q^{\pm 1}]\to 1+h\Z[[h]]$ given by mapping $aq^n\mapsto (1+nh)^a$.

If the above integrality question is asked in the sense of ``Is the
total Chern class of knots expressible as a sum of diags with integer
coefficients?'' then the answer is no, as is shown in
Section~\ref{Section:NonIntegrality}:

\begin{qthm}
  There is a degree four weight system taking integer values on all
  diagrams such that when evaluated on the total Chern class of the
  trefoil gives $-5/4$, thus the total Chern class is not
  expressible as an integer linear combination of diagrams.
\end{qthm}

This would seem to preclude the existence of a ``universal'' space
$\mathbb X$ for which the Kontsevich integral would be the Chern
character.  However, the previous result could be taken to indicate
that it might be some sort of Chern character in individual cases.  
Thinking vaguely in the realm of topological quantum field theory, the
invariant of a knot lives in the Verlinde algebra, that is
the vector space associated to a torus.  The Verlinde algebra seems to
have two manifestations, 
as the $K$-group of representations of an algebra or as the
invariant functionals on that algebra.  The isomorphism from the former
to the latter is the character map, that is the trace of the
representation.  In the case of a Lie algebra (or rather its
universal enveloping algebra) this is essentially the universal Chern
character for the corresponding group.  Dually,
the algebra $\A$ can be viewed as a ``universal'' source of
characteristic classes, giving rise to characteristic classes of
$G$-bundles for whatever $G$ and for holomorphic bundles over
holomorphic symplectic manifolds.  In any case, there is certainly
something character-esque going on.

\rk{Notation} The monoid of oriented knots equipped with the
connect sum operation is denoted $\knots$, and $\Z\knots$ denotes the
``monoid ring'' i.e.\ the ring consisting of integer linear
combinations of knots.  By $\A$ is denoted the algebra of $\Q$-power
series in connected unitrivalent diagrams modulo the STU and 1T
relations, equipped with the connect sum product.  The Kontsevich
integral $Z\colon \knots\to\A$ is normalized so that
$Z(\text{unknot})=1$, and hence is multiplicative.  Multiplicative
normalizations in general will be notated ungarnished --- $Z$, $J$, $
w_\rho$, {\em etc.} --- tildes will be used to denote {\em quantum\/}
normalizations --- $\tilde Z$, $\tilde J$, $\tilde w_\rho$, {\em etc.}
So for instance, $\tilde Z(k)=\tilde Z(\text{unknot})Z(k)$.

\rk{Acknowledgements}
Thanks go to Sir Michael Atiyah for suggesting that I read
\cite{AtiyahHirzebruch:Cohomologie-operationen}; to Jacob Mostovoy,
Justin Roberts and Elmer Rees for many conversations; and also to
Christian Kassel, Thang Le, Ted Stanford, Dylan Thurston, and the
person who refereed this by accident for useful input.  The work was
done whilst I held a Royal Society Australian Fellowship at the
University of Melbourne.

%%%%%%%%%%%%%%%%%%%%%%%%%%%%%%%%%%%%%%%%%%%%%%%%%%
%%%                                            %%%
%%%        Chern charater v. Kontsevich        %%%
%%%                                            %%%
%%%%%%%%%%%%%%%%%%%%%%%%%%%%%%%%%%%%%%%%%%%%%%%%%%

\section{The Chern character v.\ the Kontsevich integral}
\label{Section:ChernVKontsevich}

In this section I give details of the analogous properties shared by
the topological Chern character and the Kontsevich integral, which act
as the motivation for the definition of the total Chern class of knots
in Section~\ref{Section:TotalChernClass}.

%%%%%%%%%%%%%%%%%%%%%%%%
%%%  Chern Character %%%
%%%%%%%%%%%%%%%%%%%%5%%%

\subsection{The Chern character}
I will give a reminder on the topological Chern character before
giving some of its properties, more details can be found in, say,
\cite{Karoubi:KTheory}.  Suppose that $X$ is a CW complex (eg.\
a manifold, if you prefer), then $K(X)$, the $K$-group of $X$, is the
Grothendieck group of the monoid of vector bundles on $X$.  A general
element of $K(X)$ is a formal difference of two vector bundles, and
the group operation comes from $\oplus$, the Whitney sum operation on
vector bundles.  The tensor product, $\otimes$, of vector bundles
makes $K(X)$ into a ring.  The $i$th Chern class $c_i(\eta)\in
H^{2i}(X,\Z)$ of a vector bundle $\eta$ is an integral cohomology
class on $X$ which only depends on the class of $\eta$ in the
$K$-group of $X$.  The Chern character map $\ch\colon K(X)\to
H^{\text{even}}(X,\Q)$ is a special rational linear combination of
Chern classes which has the following properties.

\begin{description}
\item[Multiplicativity]  The Chern character is a ring map from a ring
to an algebra over the rationals, i.e.
$ \ch (\eta\oplus\zeta) = \ch (\eta)+\ch (\zeta)$ and
$\ch (\eta\otimes\zeta) = \ch (\eta)\cup\ch (\zeta)$.

\item[Rational associated graded map]  The $K$-group of a CW
complex is naturally filtered via the skeletal filtration of the
complex: i.e.\ if $ X^i\hookrightarrow X$ is the inclusion of the $i$th
skeleton into $X$, then the $i$th filtration, $K_i(X)$, can be defined
as $K_i(X):=\ker\{K(X)\to K(X^i)\}$.  The Chern character induces an
isomorphism on the rationalized associated graded objects:
$$ \ch\colon \bigl( K_{2i}(X)/K_{2i+1}(X)\bigr) \otimes \Q
\stackrel{\cong}{\longrightarrow} H^{2i}(X,\Q).$$
\item[Integrality]  Assuming that $H^{\blob}(X,\Z)$ is torsion-free,
$H^{\blob}(X,\Z)$ can be identified with its image under the inclusion
$H^{\blob}(X,\Z) \hookrightarrow H^{\blob}(X,\Q)$.  In this case
a rational class $x\in H^{2i}(X,\Q)$ is integral if and only if there
is an element 
$\eta\in K_{2i}(X)$ such that
$\ch(\eta)=x+\text{higher order terms}$.
\item[A result of Adams \cite{Adams:ChernCharUnitary%
%,    AtiyahHirzebruch:Cohomologie-operationen
} (briefly)]
If $x\in
H^{2i}(X,\Z)$ and  $\eta\in K_{2i}(X)$ are as above, write
$\ch(\eta)=x+\sum_{r=1}^\infty x_{2i+2r}$ with $x_{2i+2r}\in
H^{2i+2r}(X,\Q)$.  Now let  $p$ be prime
and $r$ be an integer not divisible by $p-1$ then 
$x_{2i+2r}p^{\lfloor r/(p-1)\rfloor}$ is integral
modulo $p$, and modulo $p$ it can be expressed in terms of the inverse
Steenrod reduced power of $x$, thus modulo $p$ it does not depend on
the choice of $\eta$.  
\item[Atiyah-Hirzebruch-Grothendieck-Riemann-Roch Formula]
  Supposing\break 
  that $f\colon X\to Y$ is a map between oriented manifolds one can
  define the push forward on cohomology $f_H:H^\blob(X;\Z)\to
  H^{\blob+\dim Y-\dim X}(Y;\Z)$ by using the push forward on homology
  together with Poincar\'e duality.  In fact the orientability
  assumption can be weakened so that the push-forward can be defined
  provided that the\/ {\em map\/} is orientable, ie.\ the first
  Steifel-Whitney class $w_1(f):=w_1(X)-f^*w_1(Y)\in H^1(X,Z_2)$
  vanishes.  We want to define an analogous push-forward for
  $K$-theory, this requires that the map is\/ {\em orientable in
  $K$-theory\/}, or, in other words, $\text{Spin}^c$, which means the
  following.  Let $w_2(X)\in H^2(X,\Z_2)$ denote the second
  Stiefel-Whitney class, and assume for simplicity that the dimensions
  of $X$ and $Y$ have the same parity, then the map $f$ is
  $\text{Spin}^c$ if there is an integral class $c_1(f)\in H^2(X,\Z)$ such
  that $c_1(f)\equiv w_2(X)-f^*w_2(Y) \mod 2$.  Given that $f$ is
  $\text{Spin}^c$, one can define a push forward on K-theory:
  $f_K:K(X)\to K(Y)$.  If $X$ is a smooth manifold then $\hat A(X)\in
  H^\blob(X;\Q) $ is defined in terms of the Pontryagin classes of the
  tangent bundle of $X$ and a sequence involving the power series
  $\sinh(x/2)/(x/2)$.  The Chern character naturally commutes with the
  pull-backs on homology and $K$-theory, the following theorem shows
  how it behaves with respect to the push-forward maps.  (For more
  details see, eg., \cite[Chapter V.4]{Karoubi:KTheory}.)
\end{description}

\begin{theorem}\label{AHGRR}{\rm\cite{AtiyahHirzebruch:RiemannRochDifferentiable}}\qua
  If $f\colon X\to Y$ is a map between smooth oriented manifolds whose
  dimensions have the same parity and there exists a cohomology class
  $c_1(f)\in H^2(X,\Z)$ which satisfies $ c_1(f)\equiv
  w_2(X)-f^*w_2(Y) \pmod 2$ then, for all $\eta\in K(X)$,
  $$\ch \bigl(f_K(\eta)\bigr)\cup \hat A(Y)^{-1}=f_H\bigl(\ch(\eta)\cup \hat
  A(X)^{-1}\cup \exp(c_1(f)/2)\bigr).$$
\end{theorem}

%%%%%%%%%%%%%%%%%%%%%%%%%%%%%
%%%  Kontsevich integral  %%%
%%%%%%%%%%%%%%%%%%%%5%%%%%%%%

\subsection{The Kontsevich integral}
The Kontsevich integral, $Z\colon \Z\knots\to \A$, has analogous
properties to those of the Chern character described above, as follows.
\begin{description}
\item[Multiplicativity]
  The Kontsevich integral is a ring map from a ring to an algebra over
  the rationals, i.e.\ it satisfies $Z(k+ l)=Z(k)+Z(l)$ and $Z(k\#
  l)=Z(k).Z(l)$ where $k$ and $l$ are integer linear combinations of
  knots.  In fact it is actually a bialgebra map, that is it commutes
  with the coproduct:
   $\Delta\circ Z(k) = (Z\otimes Z)\circ \Delta (k)$.

\item[Rational associated graded map]  
  The ring $\Z\knots$ is naturally graded by the Vassiliev filtration,
  this can be defined as $$(\Z\knots)_{i}=\{\text{linear combinations
  of knots with $\ge i$ double points}\},$$ knots with double points
  being considered as elements of $\Z\knots$ via the formal resolution
  $\vpic{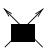}=\vpic{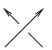}-\vpic{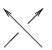}$.  A fundamental property of the
  Kontsevich integral is that it induces an isomorphism on the
  rationalized associated graded objects: 
  $$ Z\colon \bigl((\Z\knots)_{i}/(\Z\knots)_{i+1}\bigr) \otimes \Q
                            \stackrel{\cong}{\longrightarrow} \A_{i}.$$ 
  This is essentially what is meant by it being a universal Vassiliev 
  invariant.
\item[Integrality] 
 Provided that $\A_\Z$, the analogue of $\A$ defined with integer
 coefficients, is torsion-free, an element $D\in\A_i$ can be expressed
 as an {\em integer\/} linear combination of diagrams if and only if
 there is some integer linear combination, $k\in (\Z\knots)_{i}$, of
 $i$-singular knots such that $ Z(k)=D+\text{higher order terms}$.

\item[Congruences]
 I do not know if any analogue of the result of Adams holds.  However,
 the results on congruency between the terms in the Kontsevich
 integral in Proposition~\ref{congruences} seem to point in that
 direction.

\item[Cabling formula]
 The reader is directed to \cite{Willerton:KontIntAlgStructures} for
 more details; here I will sketch briefly.  For $m$ and $p$ coprime
 integers, there is the well defined notion of $(m,p)$ cabling on
 framed knots.  From a framed knot $k$, the $(m,p)$-cable $\Psi^{m,p}(k)$
 is defined so that it twists $m$ times longitudinally around the knot
 $k$ and $p$ times meridianally.  The framed Kontsevich integral, $\f
 Z\colon \Z\f\,\knots \to \f\A$, can be defined to be a map from
 integer linear combinations of framed knots to rational linear
 combinations of diagrams modulo just the STU relations.  One can
 consider $\f\A$ to be $\A$ with one generator adjoined in degree one:
 $\f\A\cong\A[[\thetadiag]]$.  The framed Kontsevich integral can be
 defined for a framed knot $k$ as $\f
 Z(k):=Z(k).\exp\left(F(k)\thetadiag/2\right)$, where $F(k)$ is the
 framing number of $k$ and $Z(k)$ is the usual Kontsevich integral of
 the underlying unframed knot.

 The wheels element, $\Omega$, of $\A$ is defined in terms of ``wheels
 diagrams'' and a ``disjoint union product'', via a power series akin
 to the square root of the $\hat A$ power series.
\end{description}
  The following theorem is the statement of how the framed Kontsevich
  integral behaves under the cabling operation.
\begin{theorem}[Le, see \cite{DThurston:Thesis}]
  If $m$ and $p$ are coprime integers and $k$ is a framed knot, then
  $$\f Z\left(\Psi^{m,p}(k)\right).\Omega = \psi^m\left(\f Z(k).\Omega.
        \exp\left(\tfrac{p}{2m}\thetadiag\right)\right).$$
\end{theorem}
\noindent The similarity with Proposition~\ref{AHGRR} is striking.

%%%%%%%%%%%%%%%%%%%%%%%%%%%%%%%%%%%%%%%%%%%%%%%%%%
%%%                                            %%%
%%%             Total Chern class              %%%
%%%                                            %%%
%%%%%%%%%%%%%%%%%%%%%%%%%%%%%%%%%%%%%%%%%%%%%%%%%%

\section{The total Chern class}
\label{Section:TotalChernClass}

In the first Section~\ref{SubSection:ClassicalChernClass} some facts
about the topological total Chern class and its relation to the Chern
character are recalled.  These are used in
Section~\ref{SubSection:KnotChernClass} to define a total Chern class
for knots.

%%%%%%%%%%%%%%%%%%%%%%%%
%%%  Classical case  %%%
%%%%%%%%%%%%%%%%%%%%5%%%

\subsection{The classical case}
\label{SubSection:ClassicalChernClass}

The Chern classes of a vector bundle over a space $X$ can be added
together to give the total Chern class $c\colon \Vect(X)\to
H^{\text{ev}}(X,\Z)$ with $c(\xi)=1+\sum_{i=1}^\infty c_i(\xi)$.  This
is multiplicative under Whitney sum of vector bundles: $c(\xi\oplus
\eta)=c(\xi)\cup c(\eta)$.  So, defining
$1+H^{\text{ev}+}(X,\Z)$ to be the group of multiplicative units
in $H^{\text{ev}}(X,\Z)$, by the universal property of the $K$-group
we get a group homomorphism
$$c\colon K(X)\to 1+H^{\text{ev}+}(X,\Z); \qquad 
  \xi\mapsto 1+\sum_{i=1}^\infty c_i(\xi), \qquad c_i(\xi)\in
            H^{2i}(X,\Z).$$
The information in the total Chern class is essentially the same as
that in the Chern character, described in
Section~\ref{Section:ChernVKontsevich}, and, modulo torsion, they can
be obtained from each other in the fashion described below.  Whereas
the Chern character is very well behaved algebraically and arises
naturally in Chern-Weil theory and index theory, it is the total Chern
class that is somehow closer to actual topology.  The Chern classes
are defined integrally and have interpretations as obstructions, in
terms of cells in Grassmanians and as being Poincar\'e dual to
submanifolds.  For instance, the zero set of a generic section of a
rank $n$ bundle $\xi$ is a codimension $2n$ submanifold which is
Poincar\'e dual to the top Chern class $c_n(\xi)\in H^{2n}(X,\Z)$.

Switching between the total Chern class (modulo torsion) and the Chern
character is entirely analogous to changing between elementary
symmetric functions and power sums, and the terms in the two classes are
related by Newton formulae (see eg.\ \cite{Karoubi:KTheory}).
  An equivalent approach is the following
formula which follows easily from the splitting principal for vector
bundles \cite{Hirzebruch:TopologicalMethods}.  If $\xi \in K(X)$ then 
$$c(\xi)\equiv \exp\bigl((-1)^{\deg/2-1}(\deg/2-1)!\, \ch (\xi)\bigr)
\mod \text{torsion},$$ where $(-1)^{\deg/2-1} (\deg/2-1)!$ is the
function on $H^{\text{ev}}(X,\Q)$ which multiplies a homogeneous
element of degree $2n$ by the factor $(-1)^{n-1} (n-1)!$, with $(-1)!$
being interpreted as 0.  Using this each Chern class (modulo torsion)
is given by a polynomial in the Chern character classes.  Eg.\
$c_2(\xi)=\tfrac 12 \ch _1(\xi)^2-\ch _2(\xi)$.

%%%%%%%%%%%%%%%%%%%%%%%%%%%%%%%
%%%  Chern class for knots  %%%
%%%%%%%%%%%%%%%%%%%%5%%%%%%%%%%

\subsection{The total Chern class for knots}
\label{SubSection:KnotChernClass}
In this section we will use the relationship between the topological
Chern character and the topological total Chern class given above,
together with the analogy between the Kontsevich integral and the
topological Chern character of Section~\ref{Section:ChernVKontsevich}
to define a ``total Chern class for knots''.  The original point of
doing this was to test the analogy and see if an ``integer'' invariant
was obtained.  As we will see in
Section~\ref{Section:LieAlgebraWeights} this is indeed the case when
a weight system coming from an irreducible representation of a simple
Lie algebra is applied. 

To define the total Chern class for knots, let $1+\A^+$ be the set of
elements in $\A$ whose term in degree zero is precisely the diagram
with no internal graph; this set forms a group under the connect sum
operation.  Now by analogy with the formula in
Section~\ref{SubSection:ClassicalChernClass}, define the total Chern
class of knots $c\colon\Z\knots\to 1+\A^+$ by
$$c(k):= \exp\bigl((-1)^{\deg-1}(\deg-1)!\, Z(k)\bigr).$$
Three remarks are in order.  The first is that the grading used here is
the conventional grading on $\A$, there are several arguments as to
why the natural grading ought to be double that, but in this paper I
will stick to the conventional grading: this is why the grading differs
by a factor of two from that in the expression for the topological
Chern class above.  The second remark is that I am ignoring any
questions of torsion here --- it is not known whether there is any
torsion in the analogue of $\A$ defined over $\Z$.  The third remark
is that if a basis of $\A$ is chosen which consists of the monomials
in some set of connected diagrams (see
\cite{Willerton:KontIntAlgStructures}), then in the transition from
the Kontsevich integral to the total Chern class, the coefficient of a
connected diagram $D$ just gets multiplied by
$(-1)^{\deg(D)-1}(\deg(D)-1)!$, but the coefficients of non-connected
diagrams get much more messed up in general.

We finish this section with two basic properties of the total Chern
class of knots.

\begin{prop}
The total Chern class of knots is a group homomorphism:
$c(k+l)=c(k).c(l)$.  
\end{prop}
\begin{proof}
This follows from the linearity of the Kontsevich integral,
$Z(k+l)=Z(k)+Z(l)$, together with the fact that the exponential map satisfies 
$\exp(A+B)=\exp(A)\exp(B)$. 
\end{proof}
Note here that the sum is the formal algebraic sum in $\Z\knots$, and
this means that $c$ is unchanged by the formal addition of trivial knots, as
$c(\text{unknot})$ is the unit in $1+\A^+$.  This is the analogue of
the fact that the topological Chern class does not detect trivial
bundles of any rank.

\begin{prop}
The total Chern class of knots
is a universal Vassiliev invariant in the following sense. 
If $k$ is a knot with $n$ double points and $D$ is the underlying
chord diagram with $n$ chords, then the total Chern class satisfies
$$c(k)=1+(-1)^{n-1}(n-1)!D+\text{higher order terms}.$$ 
\end{prop}
\begin{proof}
The Kontsevich integral is a universal Vassiliev invariant in the
sense that if $k$ satisfies the hyphotheses then $Z(k) =
D+\text{higher order terms}$.  The result follows immediately from the
definition of the total Chern class.
\end{proof}

%*********************************************
%*                                           *   
%* LIE ALGEBRAS, POLYS, AND WEIGHT SYSTEMS   *
%*                                           *   
%*********************************************

\section{Integrality on the level of Lie algebras}
\label{Section:LieAlgebraWeights}

In this section it is shown that the total Chern class of a knot gives
an integer power series when a weight system coming from an
irreducible representation of a simple Lie algebra is evaluated on
it.  Recall the notational conventions from the introduction that a
tilde means quantum normalization.

There is a standard construction in the theory of Vassiliev
invariants, by which one obtains a weight system, $\tilde w_\rho\colon
\A\to \Q[[h]]$, from the data of a representation $\rho$ of a Lie
algebra $\mathfrak g$ which is equipped with an invariant, symmetric,
bilinear form.  The idea is that from a diagram $D$ one constructs an
element $w_{\mathfrak g}(D)\in ZU(\mathfrak g)$ in the centre of the
universal enveloping algebra of the Lie algebra, one
then takes the trace in the representation $\rho$ of this element, and
multiplies by an appropriate power of $h$: i.e.\ $\tilde w_\rho(D):=
\text{Tr}_\rho w_{\mathfrak g}(D) \cdot h^{\deg D}$.

Given the same Lie algebraic data, one can form the Reshetikhin-Turaev
invariant $\tilde\tau_\rho\colon \knots\to \Q[[h]]$.  Recalling that $\tilde
Z$ is the\/ {\em quantum\/} normalization of the Kontsevich integral,
it is a fundamental theorem \cite{LeMurakami:Framed,Kassel:QuantumGroups}
of the theory of the Kontsevich integral
that the Reshetikhin-Turaev invariant $\tilde\tau_\rho$ factors as $\tilde
w_\rho\circ \tilde Z$, i.e.\ the following diagram of sets commutes.
\begin{align}\begin{array}{ccc}
\knots&\stackrel{\tilde Z}{\longrightarrow}& \A\\
&\hbox{\llap{$\tilde \tau_\rho$}}\searrow
&\downarrow\lefteqn{\tilde w_\rho}\rule[-.7em]{0em}{2em} \\
 &&\Q[[h]].
\end{array}
\end{align}

Actually, the above constructions lead to\/ {\em link} invariants and
not just knot invariants.  The modifications below will be made
because I am thinking very specifically about the case of knots.

Consider the Jones polynomial: in its quantum normalization it is a
link invariant $\tilde J\colon \text{\bf links}\to \Z[q^{\pm 1/2}]$ 
which maps the
empty link to $1$ and the unknot to $q^{1/2}+q^{-1/2}$.  
If one substitutes $q=e^h$
and expands in powers of $h$, then one obtains precisely
the Reshetikhin-Turaev invariant 
$\tilde\tau_\rho$ where $\rho$ is the fundamental representation of
$\mathfrak{su}_2$ (equipped with the bilinear form coming from the
trace in the fundamental representation).  This is fine, but if one is
thinking purely of knots then it makes sense to normalize so that the
unknot gets mapped to $1$; so to this end, set $J(k):=\tilde J(k)\tilde
J(\text{unknot})^{-1}$.  This normalization has two advantages: (i) $J(k)$
is a genuine Laurent polynomial in $q$ with no odd powers of $q^{1/2}$;
and (ii) it is multiplicative
under connected sum, i.e.\ $J(k\# l)=J(k)J(l)$.

In general, for $\rho$ an irreducible representation of a simple
Lie algebra, one has the quantum link invariant $\tilde J_\rho\colon
\text{\bf links}\to \Z[q^{\pm 1/M}]$, which takes values in the ring of
Laurent polynomials in some fractional power, $q^{1/M}$, of $q$, where
$M$ is some integer depending on the representation.  This is related
to the Reshetikhin-Turaev invariants via the following commutative diagram:
\begin{align}\begin{array}{ccc}
\text{\bf links}&&\\
\llap{{$\tilde J_\rho$}}\downarrow\rule[-.7em]{0em}{2em}
&\searrow\lefteqn{\tilde \tau_\rho}&\\
 \Z[q^{\pm 1/M}]&\stackrel{q\mapsto e^h}{\longrightarrow}&\Q[[h]].
\end{array}
\end{align}

 If attention is
restricted to knots and the invariant $\tilde J_\rho$ is renormalized to
$J_\rho$ so that the
unknot is mapped to $1$ then one can ask what the image of the map $J_\rho$
is.  One might expect that, as in the case of the Jones polynomial
above, the image is contained in $\Z[q^{\pm 1}]$.  However, I was
surprised to find nothing in the literature on this question, until
the recent paper of Le \cite{Le:IntegralityAndSymmetry}.  Perhaps
one reason for this lack is that the question is not terribly natural
when the approach to the quantum invariants is skein theoretic,
tangle theoretic, or via Markov traces.  In the 
current context, the question is very natural.  The following theorem
is proved by a thorough understanding of Lusztig's work on quantum groups.

\begin{theorem}[Le \cite{Le:IntegralityAndSymmetry}] If $\rho$ is an
  irreducible representations of a simple Lie algebra, and $J_\rho$ is
  normalized to take value $1$ on the unknot, then it gives a
  multiplicative map $J_\rho\colon \knots\to \Z[q^{\pm 1}]$.
\end{theorem}

For $\rho$ an irreducible representations of simple Lie algebra
$\mathfrak g$, define a weight system in the following manner.  As the
universal $\mathfrak g$ weight system $w_{\mathfrak g}$ takes values
in the centre of the universal enveloping algebra and the latter acts
by scalars in all irreducible representations, define $w_\rho$ to be
precisely the value of this scalar.  As $\tilde w_\rho$ was defined as
the trace of this scalar operator, it is clear that $\tilde
w_\rho=(\dim\rho) .w_\rho$.  This normalization is algebraically
well-behaved in the sense of the following theorem.
\begin{theorem}  
  If $\rho$ is an irreducible representation of a Lie algebra then 
  the weight system $w_\rho\colon \A\to \Q[[h]]$ is an algebra map.
\end{theorem}
\begin{proof}
  This follows immediately from the fact that $w_{\mathfrak g}$ is
  multiplicative, and the fact that $\rho$ is multiplicative on the
  universal enveloping algebra.
\end{proof}
Now extend $J_\rho$ linearly to $\Z\knots$, and define the ring map
$\ch\colon \Z[q^{\pm 1}]\to \Q[[h]]$ by $q\mapsto e^h$; the notation
will be explained after the following theorem.
\begin{theorem}
  If $\rho$ is an irreducible representation of a simple Lie algebra,
  then the following is a commutative diagram of filtered rings:
  \begin{align*}\begin{array}{ccc}
    \Z\knots&\stackrel{Z}{\longrightarrow}& \A\\
     \llap{{$J_\rho$}}\downarrow&&
        \downarrow\lefteqn{w_\rho}\rule[-.7em]{0em}{2em}\\ 
     \Z[q^{\pm 1}] & \stackrel{\ch}{\longrightarrow} &\Q[[h]].
     \end{array}
  \end{align*}
\end{theorem}
\begin{proof}
  Let $k$ be a knot.  By the commutativity of diagrams (1) and (2)
  above, for the quantum normalizations, $\ch\bigl( \tilde
  J_\rho(k)\bigr)=\tilde w_\rho\bigl( \tilde Z(k) \bigr)$.  This can
  be written in the multiplicative normalizations as
\begin{align*}
\ch \bigl(\tilde J_\rho(\text{unknot}).J_\rho(k)\bigr)
   &=(\dim\rho). w_\rho \bigl( \tilde Z(\text{unknot}).Z(k)\bigr).
\intertext{By multiplicativity of $\ch$ and $ w_\rho$,}
\ch \bigl(\tilde J_\rho(\text{unknot})\bigr).\ch \bigl(J_\rho(k)\bigr)
   &=(\dim\rho). w_\rho \bigl( \tilde Z(\text{unknot})\bigr)
                                 .w_\rho \bigl(Z(k)\bigr)\\
&=\tilde  w_\rho \bigl( \tilde Z(\text{unknot})\bigr)
                                 .w_\rho \bigl(Z(k)\bigr)
\end{align*}
But $\ch \bigl(\tilde J_\rho(\text{unknot})\bigr)= \tilde w_\rho
\bigl( \tilde Z(\text{unknot})\bigr)$ and this is an invertible
element in $\Q[[h]]$, so $\ch\bigl( J_\rho(k)\bigr)=w_\rho\bigl(
Z(k)\bigr)$ as required.
\end{proof}

As $Z$ is being thought of as some kind of Chern character, one can
ask what my interpretation of $\ch\colon\Z[q^{\pm 1}] \to \Q[[h]]$,
$q\mapsto e^h$, is in this language.  Thinking na\"\i vely as a
topologist, $\Z[q^{\pm 1}]$ can be thought of as the Grothendieck
group of vector bundles over infinite projective space $\mathbb C P^\infty$,
with $q$ representing the canonical line bundle; similarly, $ \Q[[h]]$
can be thought of as the (completed) rational, ordinary cohomology of
$\mathbb C P^\infty$, with $h$ being the canonical generator in degree
$2$ which is the first Chern class of the canonical line bundle $q$.
In this framework $q\mapsto e^h$ corresponds\/ {\em precisely\/} to
the Chern character --- this is the reason that I called it $\ch$.

The natural question to ask next is: what corresponds to the total
Chern class?  This is easy.  Let $1+h\Z[[h]]$ be the group of
multiplicative units of the ring of formal power series $\Z[[h]]$ and
define $c\colon \Z[q^{\pm 1}] \to 1+h\Z[[h]]$ to be the map of\/ {\em
  groups\/} given by $q^n\mapsto 1+nh$, i.e.\ $c(a+b)=c(a)c(b)$ for
all $a,b \in\Z[q^{\pm 1}]$, so for instance, $2q^{-3}-q\mapsto
(1-3q)^2(1-h+h^2-h^3+\dots)$. 
\begin{theorem}
  For $\rho$ an irreducible represention of a simple Lie algebra, the
  following is a commutative diagram of groups:
  \begin{align*}\begin{array}{ccc}
  \Z\knots&\stackrel{c}{\longrightarrow}&1+ \A^+\\
  \llap{{$J_\rho$}}\downarrow&&
        \downarrow\lefteqn{w_\rho}\rule[-.7em]{0em}{2em}\\ 
  \Z[q^{\pm 1}] & \stackrel{c}{\longrightarrow}\ 
   1+h\Z[[h]]\ \hookrightarrow
   &1+h\Q[[h]].
   \end{array}
  \end{align*}
\end{theorem}
\begin{proof}This follows from the multiplicativity of the weight
  system $w_\rho$ and the fact that both the maps labelled $c$ can be
  expressed in terms, respectively, of the maps $Z$ and $\ch$ via the
  same formula which is polynomial in each degree.
\end{proof}
The following corollary is immediate and is the expression of the
integrality of the total Chern class for knots on the level of Lie
algebras. 
\begin{cor}\label{LieIntegrality}
  The image of the total Chern class of knots is integral as far as
  irreducible representations of simple Lie algebras can detect, i.e.\ 
  if $\rho$ is such a representation, then $w_\rho\circ
  c(\Z\knots)\subset \Z[[h]]$.
\end{cor}

It is perhaps worth contemplating, in view of the above, whether the
quantum invariants, such as the Jones polynomial, have any K-theoretic
interpretation.  Note that, thought of as a torsion, the Alexander
polynomial takes values in something akin to a Whitehead group.

%*********************************************
%*                                           *   
%*        LOW ORDER CALCULATIONS             *
%*                                           *   
%*********************************************

\section{Non-integrality on the level of diagrams}
\label{Section:NonIntegrality}

It will be shown here that the Chern class does not live in the
integer lattice of $\A$, despite the integrality results of the
previous section.

\begin{theorem}\label{DiagramNonIntegrality}
  There is a degree four weight system taking integer values on all
  diagrams such that when evaluated on the total Chern class of the
  trefoil gives $-5/4$, thus the total Chern class is not
  expressible as an integer linear combination of diagrams.
\end{theorem}
\begin{proof}
In the appendix an expression for the Kontsevich integral is given in
terms of integer
valued knot invariants, $B$, $C$, $D_1$,  $D_2$, $E_1$, $E_2$ and
$E_3$.
If the transformation to the total
Chern class is made, then one finds:
\begin{align*}
c(k)&=\vpic{O}-B(k) \vpic{tI} + 2C(k) \vpic{t2I} \\
    &\qquad -\tfrac{1}{4}\left(D_1(k) \vpic{t3I} + 
             D_2(k) \vpic{w4} + 10B(k)^2 \vpic{tItI}\right)\\
    &\qquad + 2E_1(k) \vpic{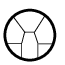} + E_2(k) \vpic{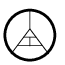}+ E_3(k)
             \vpic{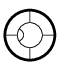}+22B(k)C(k)\vpic{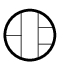}+\dots .
\end{align*}
The problem here is with the degree four piece.  One can check that the
degree four weight system that maps
$$\vpic{t3I}\mapsto 0;\quad \vpic{w4}\mapsto 1;\quad
\vpic{tItI}\mapsto 0,$$
is integer valued on all degree four diagrams.  Evaluating this on $c(k)$
one obtains the number $-D_2(k)/4$.  From the example of the trefoil in
the introduction, it is seen that $D_2(\text{trefoil})=5$ and hence
the total Chern class of the trefoil is not in the integer lattice of
$\A$.  
\end{proof}
\noindent Note that there could be things that stop this from being a
``nice'' weight system; for example I do not know if this weight
system extends to an integral multiplicative weight system, perhaps
some power of the wheel with four legs is a multiple of some integer
combination of other diagrams.

 It is instructive here
to consider the evaluation of the total Chern class under the weight
systems coming from the 
irreducible representations of $\mathfrak{su}_2$; from the previous
section this is known to give an integer power series. 
In the $d$-dimensional irreducible representation of
$\mathfrak{su}_2$, the Casimir acts by $\lambda:=(d^2-1)/2$, and, by
\cite{ChmutovVarchenko}, the corresponding weight system, $w_d$, behaves as
follows: 
%
% For $\lambda$ a
%half-integer, let $w_\lambda$ be the weight system coming from the
%representation labeled by $ \lambda$.  This gives \cite{ChmutovVarchenko}
%
$$\vpic{t3I}\mapsto 8\lambda h^4;\quad \vpic{w4}\mapsto 8\lambda^2h^4;
\quad \vpic{tItI}\mapsto 4\lambda^2h^4.$$
Then the coefficient of
$h^4$ in $w_d c(k)$ is given by
$-2\lambda((5B(k)^2+D_2(k))\lambda+D_1(k))$.  This is integral when
$\lambda$ is half-integral precisely because of the congruence
$B(k)\equiv D_2(k) \pmod 2$ from Proposition~\ref{congruences} in
the appendix.  Note
that the integrality result for the Alexander-Conway weight system
corresponds to the same congruence.

%*********************************************
%*                                           *   
%*             APPENDIX                      *
%*                                           *   
%*********************************************

%%%%%%%%%%%%%%%%%%%%%%%%%%%%%%%%%%%%%%%%%%%%%%%%%%%%%%%%%%%%
%%%                                                      %%%
%%%                 CALCULATIONS                         %%%
%%%                                                      %%%
%%%%%%%%%%%%%%%%%%%%%%%%%%%%%%%%%%%%%%%%%%%%%%%%%%%%%%%%%%%%

\section*{Appendix: The Kontsevich integral up to degree five}
In this appendix an expression for the Kontsevich integral up to degree
five is given in terms of integral knot invariants.  A corollary of
the proof is that these invariants satisfy various congruences.

\begin{theorem}
  There exist integer
valued knot invariants, $B$, $C$, $D_1$,  $D_2$, $E_1$, $E_2$ and
$E_3$ such that 
%Diagrams saved at 30% normal size
\begin{align*}Z(k)&=\exp\biggl(B(k) \vpic{tI} + C(k) \vpic{t2I}
    +\tfrac{1}{24}\left( D_1(k) \vpic{t3I} + D_2(k) \vpic{w4}
    \right)\\&\qquad\qquad  
   +\tfrac{1}{24}\left(2E_1(k) \vpic{t4I} + E_2(k) \vpic{x3tI}+ E_3(k)
   \vpic{tw4} 
    \right)+ \dots \biggr).
\end{align*}
\end{theorem}
\begin{proof}
  It follows from \cite{Willerton:KontIntAlgStructures} that the
  Kontsevich integral has this form where the invariants are rational
  valued.  Stanford \cite{Stanford:Programs} has calculated a basis
  for rational valued additive Vassiliev invariants up to degree six
  which consists of {\em integer}\/ valued invariants.  To express the
  above invariants in terms of Stanford's invariants it suffices to
  calculate the Kontsevich integral for four suitably chosen knots,
  compare the values with Stanford's table and then solve the
  requisite linear equations.  I calculated the Kontsevich integral up
  to degree five of the knots $3_1$, $5_1$, $5_2$, and $7_2$ by using
  my formulae for torus knots and Bar-Natan's Mathematica program (see
  \cite{Bar:NAT}).  Denote Stanford's invariants $II$, $III$,
  $IV_1$, $IV_2$, $V_1$, $V_2$ and $V_3$.  Define an alternative
  integral basis in degree five by $W_1=(V_1-2V_2)/9$, $W_2=V_2$ and
  $W_3=(V_2-V_3)/9$.  Then one finds

  { \footnotesize
  \begin{align*}
    B&=-II, & II&=-B,\\ C&=III, & III&=C,\\ D_1&=17II-48IV_1-24IV_2, &
    IV_1&=\tfrac{1}{12}(D_2-7B),\\ D_2&=-7II+12IV_1, &
    IV_2&=\tfrac{1}{24}(11B-D_1-4D_2),\\
     E_1&=11III-140W_1-76W_2+340W_3, &
    W_1&=\tfrac{1}{24}(22C+2E_1+7E_2+3E_3),\\
    E_2&=III-8W_1 -4W_2+16W_3,&
    W_2&=\tfrac{1}{12}(-23C-2E_1-40E_2-5E_3),\\
    E_3&=-17III+120W_1+60W_2-264W_3,&
    W_3&=\tfrac{1}{24}(-2C-15E_2-E_3).
  \end{align*}
  }
  \rm The left-hand set of equations together with the integrality of
  Stanford's invariants proves the theorem.
\end{proof}

It is worth remarking that it is natural to consider the denominators
of the Kontsevich integral, where an element $D\in\A$ has\/ {\it
  denominator\/} dividing $m\in\Z$ if $mD$ can be expressed as an
integer linear combination of diagrams.  Le \cite{Le:Denominators} has
considered this for the quantum normalization, $\tilde Z$, and found
that the denominator of the degree $n$ piece divides
$(2!3!...n!)^4(n+1)!$.  I suspect that the denominators in the
multiplicative normalization are better behaved.  Note that if the
degree $n$ piece of $\ln Z$ has denominator dividing $n!$, then the
degree $n$ piece of $Z$ has denominator dividing $n!$.  From the above
theorem, one might therefore be tempted to conjecture precisely that.

As a corollary of the above proof, one also obtains the following
intriguing result.
\begin{prop}
  \label{congruences}
  The integral knot invariants $B$, $C$, $D_1$, $D_2$, $E_1$, $E_2$
  and $E_3$ satisfy the following congruences:
  \begin{align*}
    B&\equiv D_1 \equiv D_2 \pmod{3};&
    D_1&\equiv -B \pmod 8;\\
    4C&\equiv D_2-D_1 \pmod{8};&
    D_2&\equiv -B \pmod 4;\\
    E_1&\equiv E_2+E_3 \pmod{3};&
    E_2&\equiv-E_3\pmod 8;\\
    C&\equiv E_3 \pmod{3};& C&\equiv -E_1\equiv E_2 \pmod
    4.
  \end{align*}
\end{prop}
\begin{proof}
  These congruences are consequences of the right-hand set of equations
  in the previous proof together with the fact that all of the
  invariants are integer valued.  The third congruence also requires
  Stanford's observation that $IV_1-III$\/ is always even.
\end{proof}
From these one can obtain for instance, $D_1-D_2\equiv 0 \pmod{12}$ and
$C\equiv (D_1-D_2)/4 \pmod 2$; but I have no idea what these mean.
Perhaps it is useful to note that $B \pmod 2$ is the Arf invariant.

%*********************************************
%*                                           *   
%*           BIBLIOGRAPHY                    *
%*                                           *   
%*********************************************

\Addresses\recd


\newpage

\begin{thebibliography}

\bibitem{Adams:ChernCharUnitary}
\textbf{J\,F Adams}, \emph{On {C}hern characters and the structure of the
  unitary group}, Math. Proc. Cambridge Philos. Soc. 57 (1961) 189--199

\bibitem{AtiyahHirzebruch:RiemannRochDifferentiable}
\textbf{M\,F Atiyah}, \textbf{F Hirzebruch}, \emph{Riemann-{R}och theorems for
  differentiable manifolds}, Bull. Amer. Math. Soc. 65 (1959) 276--281

\bibitem{AtiyahHirzebruch:Cohomologie-operationen}
\textbf{M\,F Atiyah}, \textbf{F Hirzebruch}, \emph{Cohomologie-operationen und
  {C}harakteristische {K}lassen}, Math. Z. 77 (1961) 149--187, in German

\bibitem{Bar:OnVIs}
\textbf{D Bar-Natan}, \emph{On the {V}assiliev knot invariants}, Topology 34
  (1995) 423--472

\bibitem{Bar:NAT}
\textbf{D Bar-Natan}, \emph{Non-associative tangles}, from: ``Geometric
  topology (Athens, GA, 1993)'', Amer. Math. Soc., Providence, RI (1997)
  139--183

\bibitem{ChmutovVarchenko}
\textbf{S\,V Chmutov}, \textbf{A\,N Varchenko}, \emph{Remarks on the
  {V}assiliev knot invariants coming from ${\rm sl}_2$}, Topology 36 (1997)
  153--178

\bibitem{Hirzebruch:TopologicalMethods}
\textbf{F Hirzebruch}, \emph{Topological methods in algebraic geometry},
  Springer-Verlag, Berlin (1995)

\bibitem{Karoubi:KTheory}
\textbf{M Karoubi}, \emph{${K}$-theory}, Grundlehren der Mathematischen
  Wissenschaften 226, Springer-Verlag, Berlin (1978)

\bibitem{Kassel:QuantumGroups}
\textbf{C Kassel}, \emph{Quantum groups}, Graduate Texts in Mathematics 155,
  Springer-Verlag, New York (1995)

\bibitem{Le:Denominators}
\textbf{T\,Q\,T Le}, \emph{On denominators of the {K}ontsevich integral and the
  universal perturbative invariant of $3$-manifolds}, Invent. Math. 135 (1999)
  689--722

\bibitem{Le:IntegralityAndSymmetry}
\textbf{T\,Q\,T Le}, \emph{Integrality and symmetry of quantum link
  invariants}, Duke Math. J. 102 (2000) 273--306

\bibitem{LeMurakami:Framed}
\textbf{T\,Q\,T Le}, \textbf{J Murakami}, \emph{The universal
  {V}assiliev-{K}ontsevich invariant for framed oriented links}, Compositio
  Math. 102 (1996) 41--64

\bibitem{Stanford:Programs}
\textbf{T Stanford}, \emph{Computer programs and data tables}, available via\nl
  {\tt ftp://geom.umn.edu/pub/contrib/vassiliev.tar.Z} (1992)

\bibitem{DThurston:Thesis}
\textbf{D\,P Thurston}, \emph{Wheeling: A diagrammatic analogue of the {D}uflo
  isomorphism}, Ph.D. thesis, U.C. Berkeley (2000), {\tt math.QA/0006083}

\bibitem{Willerton:KontIntAlgStructures}
\textbf{S Willerton}, \emph{The {K}ontsevich integral and algebraic structures
  on the space of diagrams}, from: ``Knots in Hellas '98'', Series on Knots and
  Everything 24, World Scientific (2000)  530--546

\end{thebibliography}
\end{document}